\documentclass[12pt]{article}
\usepackage{amsmath}
\usepackage{amssymb}

\begin{document}
\title{A Remark on the Kelvin Transform for a Quasilinear Equation}
\author{Peter Lindqvist}
\date{Department of Mathematics\\Norwegian University of Science and Technology}
\maketitle

The p-harmonic equation $\mathrm{div}(|\nabla u|^{p-2} \nabla u) =
0$ has certain applications to quasiregular mappings [B-I], to
non-linear potential theory, where the corresponding variational
integral
defines the $p$-Capacity, and in Physics it describes certain
non-Newtonian fluids. The equation is invariant under Euclidean motions
and the aim of this note is to study the behaviour under reflections
in spheres.

The Euler--Lagrange equation for the variational integral
\begin{equation}
	\int \, |\nabla u(x)|^pdx, \qquad
	(1 < p < \infty)
\end{equation}
$dx$ denoting the n-dimensional Lebesgue measure, is the p-harmonic
equation
\begin{equation}
	\mathrm{div}(|\nabla u|^{p-2}\nabla u) = 0 \qquad 
	(1 < p < \infty).
\end{equation}
The solutions are called p-{\it harmonic} {\it functions}. Usually
Eqn (2) is interpreted in the weak sense: a function $u$ in the
Sobolev space $W^{1,p}_{\mathrm{loc}}(G), \, G$ being a domain in the
Euclidean n-dimensional space $\mathbb{R}^n$, is p-harmonic in $G$, if
\begin{equation}
	\int \,|\nabla u|^{p-2}\nabla u\cdot\nabla\eta\, dx = 0
\end{equation}
whenever $\eta \in C^{\infty}_0(G)$. It is known that the weak solutions
are (equivalent to functions) of class $C^{1,\alpha}_{\mathrm{loc}}(G)$.
This means that the gradient $\nabla u$ is locally Holder continuous
with the exponent $\alpha = \alpha(n,p) > 0$, cf. [D,E,Lew,T,U].
(In the
plane the best Holder exponent is known [I-M].)

The integral $\int |\nabla u|^p$ is conformally invariant when
$p = n$ ("the border-line case"). This important observation was
made by Ch. Loewner in 1959, cf. [Loe]. In the n-dimensional space
this means that the n-harmonic functions are invariant under Mobius
transformations. See [B, Chapter 3] and [A] for Mobius transformations
in $\mathbb{R}^n$. Since all p-harmonic functions are preserved under
Euclidean motions and homothetic transformations, the crucial point
is easily reduced to the reflection (or inversion)
	$$
	x^{\ast} = \frac{x }{|x|^2}
	$$
in the unit sphere and the proof is a matter of a simple calculation.
If $u$ is n-harmonic in the domain $G$ in $\mathbb{R}^n$, so is
	\begin{equation}
	v = v(x) = u\Bigl(\frac{x}{| x|^2}\Bigr)
	\end{equation}
in the reflected domain $G^{\ast}$. If necessary, define $0^{\ast} =
\infty$.

If $p \not= n$, p-harmonicity is not, in general, preserved under
reflections in spheres: the Mobius invariance brakes down.

For $p = 2$ Eqn (2) reduces to the Laplace equation $\Delta u = 0$.
Even though harmonic functions in higher dimensions are not preserved
under reflections in spheres, this deficiency is fortunately
compensated by the celebrated Kelvin transform\footnote{W.
Thomson (Lord Kelvin), Journal de math\'ematiques pures et
appliqu\'ees $\underline {12}$ (1847), p. 256.}: if $u$ is harmonic
in
$G$, then
	\begin{equation}
	v = v(x) = |x|^{2-n} u \Bigl(\frac{x}{|x|^2}\Bigr)
	\end{equation}
is harmonic in the reflected domain $G^{\ast}$. See [H, Section 1.9].
A recent application is given in [Leu].

Note that formulae (4) and (5) can be written as
	\begin{equation}
	v = v(x) = | x|^{\frac{p-n }{ p-1}} u
	\Bigl( \frac{x }{|x|^2}\Bigr)
	\end{equation}
valid for $p = 2$ and $p = n$. One might ask whether there is, for
general $p$, some counterpart to the Kelvin transform. Although
there is intuitive support in favour of a positive conjecture, the
answer is plainly ''no''. To some extent this answer is unexpected,
to say the least. This constitutes a serious obstacle for the
development of some parts of the theory for p-harmonic functions.

\bigskip

{\bf THEOREM} {\sl For a given $p, 1 < p < \infty$, there is
no radial function $\rho \not\equiv 0$ such that
	\begin{equation}
	v = v(x) = \rho(|x|)\,u\Bigl(\frac {x }{| x|^2}\Bigr)
	\end{equation}
is p-harmonic in the reflected domain $G^{\ast}$, whenever $u$ is
p-harmonic in $G, G \subset \mathbb{R}^n$, except in the cases $p = 2$ and
$p
= n$.}

\bigskip
{\bf Proof.} Let us first point out that a direct calculation of
$\mathrm{ div}(|\nabla v|^{p-2} \nabla v)$ does not lead to anything
comprehensible. A more ''experimental approach'' yields the
non-existence of $\rho$: transform sufficiently many p-harmonic
functions according to (7) to get a contradiction.

Transforming the p-harmonic function
	\begin{equation}
	u(x) = \begin{cases} 
	| x|^{(p-n)/ (p-1)}, \quad \text{if} \quad p \not= n, \\
	\ln(|x|),\quad \text{if} \quad p = n,
\end{cases}
	\end{equation}
we immediately get the necessary condition
	\begin{equation}
	\rho(|x |) = | x|^{(p-n)/ (p-1)},
	\end{equation}
since, as it is easily seen, the only radial p-harmonic functions
are essentially those in (7), variants like $A|x-a|^{(p-n)/
(p-1)} + B$ being included. Thus we are back at (6), indeed.

A routine calculation shows that the functions
	$$
	u = u(x_1,...,x_n) =
	(x^2_1+x^2_2+\cdots +x^2_j)^{(p-j)/2(p-1)},
	$$
$j = 1,2,...,n$ are p-harmonic, when $x^2_1 + x^2_2 +\cdots
+x^2_j \not= 0$. But transforming them according to (6), we shall
obtain p-harmonic functions only for $p = n,\, p = 2$, and $j = n$.
Indeed, it is here sufficient, although not virtually simpler, to
transform only the first function $u(x) = x_1$ in order to arrive at
a contradiction with the necessary choice (9). Now the following
technical lemma concludes our proof.

\bigskip
{\bf LEMMA} {\sl The function $v = v(x) = |x|^{\alpha} x_1$ is
p-harmonic only in the following cases: (i) $\alpha = 0$, (ii) $\alpha
= -2$ and $p = n$, (iii) $\alpha = -n$ and $p = 2$, and (iv) $\alpha
= -1$ and $p = 3-n$.}

\bigskip

{\bf Proof.} A lengthy routine calculation gives the expression
	\begin{gather*}
	\mathrm{ div}(|\nabla v|^{p-2}\nabla v) =
	|\nabla v|^{p-4}\alpha |x|^{3\alpha -4} x_1
	\lbrace [\alpha + n + (p-2)(2\alpha +3)]|x|^2 \\
	+ (\alpha + 2)[\alpha(\alpha +n) + (p-2)
	(\alpha^2 +\alpha - 1)]x^2_1\rbrace ,
	\end{gather*}
where
	$$
	|\nabla v|^2 = |x|^{2(\alpha -1)}
	\lbrace \alpha(\alpha + 2)x^2_1 + |
x|^2\rbrace.
	$$
This yields the desired result. (Case (iv) does not occur for $1
< p <  \infty$.) 

\bigskip

{\bf Epilogue.} Of course one might try to replace the reflection by
some kind of a distorted inversion like
	$$
	x^{\ast} = \frac{x}{ |x|^\beta}
	$$
where $\beta = \beta (n,p) > 0$. Indeed, the disappointing news
is that no ''Kelvin transform'' even of the type
	$$
	v = v(x) = \rho( |x|)\, u\Bigl(\frac{x}{ |x|^\beta}
	\Bigr)
	$$
exists for all p-harmonic functions $u$ if $p \not= n$ and $p \not=
2$. The proof is similar to the above one.

This leaves us with little or no hope to map an unbounded domain
onto a bounded one in the p-harmonic setting. However, sometimes one
can proceed as in [L1], using p-superharmonic functions. In the complex
plane there is, fortunately, a rich structure partly compensating
for this lack [L2].

\bigskip

{\bf References}
\begin{description}
\item{[A]}
AHLFORS, L.: Mobius Transformations in Several Dimensions, University
of Minnesota, Minnesota 1981.

\item{[B]}
BEARDON, A.: The geometry of Discrete Groups, Springer-Verlag, New
York, 1983.

\item{[B-I]}
BOJARSKI, B \& IWANIEC, T.: Analytical foundations of the theory of
quasiconformal mappings in {\bf R}$^n$, Annales Academiae Scientiarum
Fennicae Ser. A. I. Math. $\underline 8$ (1983), pp. 257-324. 

\item{[D]}
DIBENEDETTO, E.: $C^{1+\alpha}$ local regularity of weak solutions
of degenerate elliptic equations, Nonlinear Analysis TMA $\underline
7$ (1983), pp. 827-850.

\item{[E]}
EVANS, L.: A new proof of local $C^{1,\alpha}$ regularity for solutions
of certain degenerate elliptic P.D.E., Journal of Differential Equations
$\underline {45}$ (1982), pp. 356-373.

\item{[H]}
HUA, L-K.: Starting with the Unit Circle, Springer-Verlag, New York
1981.

\item{[I-M]}
IWANIEC, T. \& MANFREDI, J.: Optimum regularity for p-harmonic functions
on the plane (Manuscript).

\item{[Leu]}
LEUTWILER, H.: On a distance invariant under Mobius transformations
in {\bf R}$^n$, Annales Academiae Scientiarum Fennicae Ser. A. I.
Math. $\underline {12}$ (1987), pp. 3-17.

\item{[Lew]} LEWIS, J.: Regularity of the derivatives of solutions
to certain degenerate elliptic equations, Indiana Univ. Math. J.
$\underline {32}$ (1983), pp. 849-858.

\item{[L1]}
LINDQVIST, P.: On the growth of the solutions of the differential
equation $\mathrm{div}(|\nabla u|^{p-2}\nabla u) = 0$ in n-dimensional
space, Journal of Differential Equations $\underline {58}$ (1985),
pp. 307-317.

\item{[L2]}
LINDQVIST, P.: On p-harmonic functions in the complex plane and
curvature, Israel Journal of Mathematics $\underline {63}$ (1988),
pp. 257-269.

\item{[Loe]}
LOEWNER, C.: On the conformal capacity in space, Journal of Mathematics
and Mechanics $\underline 8$ (1959), pp. 411-414.

\item{[T]}
TOLKSDORF, P.: Regularity for a more general class of quasilinear
elliptic equations, Journal of Differential Equations $\underline{51}$
(1984), pp. 126-150.

\item{[U]}
URAL'TSEVA, N.: Degenerate quasilinear elliptic systems, Zap. Nau\v
cn. Sem. Leningrad. Otdel. Mat. Inst. Steklov $\underline 7$ (1968),
pp. 184-222. (In Russian)
\end{description}
\end{document}